%% file: main.tex
\documentclass[12pt,a4paper,leqno]{article}

\usepackage{amssymb,exscale}
\usepackage[centertags]{amsmath}
\usepackage{graphicx}
\usepackage{theorem}

\input{preamble}

\begin{document}

\title{Scaling limit of Fourier-Walsh coefficients\\(a framework)}
\author{Boris Tsirelson}
\date{}

\maketitle

\begin{abstract}
Independent random signs can govern various discrete models that
converge to non-isomorphic continuous limits. Convergence of
Fourier-Walsh spectra is established under appropriate conditions.
\end{abstract}

\input{intro}

\input{first}

\input{second}

\input{third}

\input{fourth}

\bigskip
\filbreak
\begingroup
{
\small
\begin{sc}
\parindent=0pt\baselineskip=12pt
\def\emailwww#1#2{\par\qquad {\tt #1}\par\qquad {\tt #2}\medskip}

School of Mathematics, Tel Aviv Univ., Tel Aviv
69978, Israel
\emailwww{tsirel@math.tau.ac.il}
{http://math.tau.ac.il/$\sim$tsirel/}
\end{sc}
}
\filbreak

\endgroup

\end{document}

%% file: preamble.tex
\numberwithin{equation}{section}

\theoremstyle{change}
\theorembodyfont{\rmfamily}
\newtheorem{theorem}[equation]{Theorem}
\newtheorem{lemma}[equation]{Lemma}

\newtheorem{Def}[equation]{Definition}
\newtheorem{example}[equation]{Example}

\newtheorem{note}[equation]{Note}

\renewcommand{\phi}{\varphi}

\renewcommand{\(}{\bigl(}
\renewcommand{\)}{\bigr)\vphantom{)}}

\newcommand{\Om}{\Omega}

\newcommand{\eps}{\varepsilon}

\newcommand{\de}{\delta}

\newcommand{\F}{\mathcal F}

\newcommand{\la}{\lambda}

\newcommand{\E}{\mathbb E}

\renewcommand{\Pr}[1]{\,\mathbb P\,\(\,#1\,\)\,}
\newcommand{\R}{\mathbb R}

\newcommand{\Z}{\mathbb Z}

\newcommand{\cE}[2]{\mathbb{E}\,\(\,#1\,\big|\,#2\,\)\,}

\newcommand{\qed}{\hfill$\square$}

\newcommand{\sif}{$\sigma$-field}

\newcommand{\Sra}{S^{\text{random}}}
\newcommand{\rhoKR}{\rho_\mathrm{KR}}

%% file: intro.tex
\section*{Introduction}

\begin{figure}[b!]
\begin{center}
\setlength{\unitlength}{1cm}
\begin{picture}(12,5)
\put(-0.3,3.4){\includegraphics{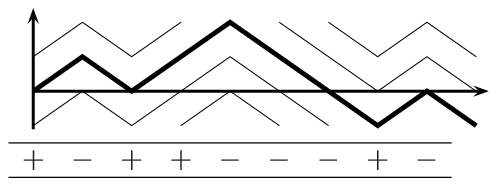}}
\put(2,3.25){\makebox(0,0){(a)}}
\put(5.65,3.4){\includegraphics{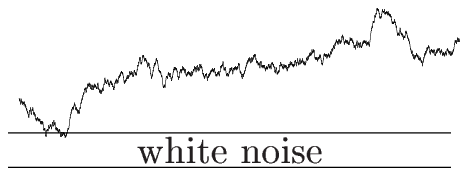}}
\put(8,3.25){\makebox(0,0){(b)}}
\put(-0.7,-0.1){\includegraphics{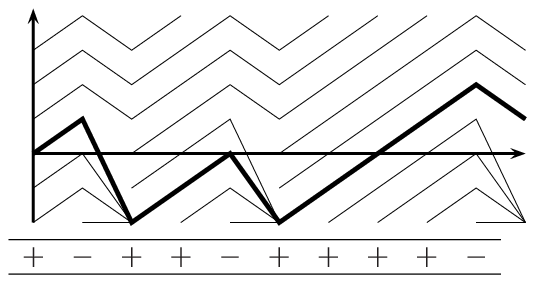}}
\put(2,-0.25){\makebox(0,0){(c)}}
\put(5.5,-0.1){\includegraphics{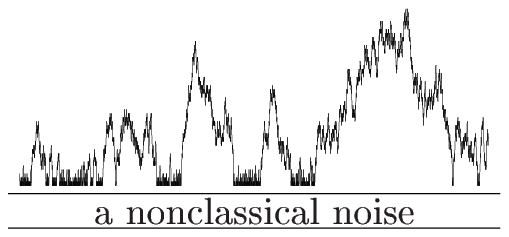}}
\put(8,-0.25){\makebox(0,0){(d)}}
\end{picture}
\caption[]{\label{fig1}\small
One-dimensional array of random signs can govern a random
walk (a) that converges in distribution to Brownian motion (b).
Such array can also govern a sticky walk (c) converging to sticky
Brownian motion (d).}
\end{center}
\end{figure}

\begin{figure}[htbp]
\begin{center}
\setlength{\unitlength}{1cm}
\begin{picture}(10,2.6)
\put(-0.5,0.3){\includegraphics{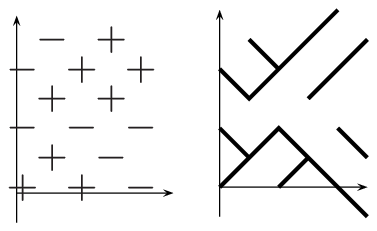}}
\put(1.75,0){\makebox(0,0){(a)}}
\put(6,0.3){\includegraphics[width=5cm,height=2.25cm]{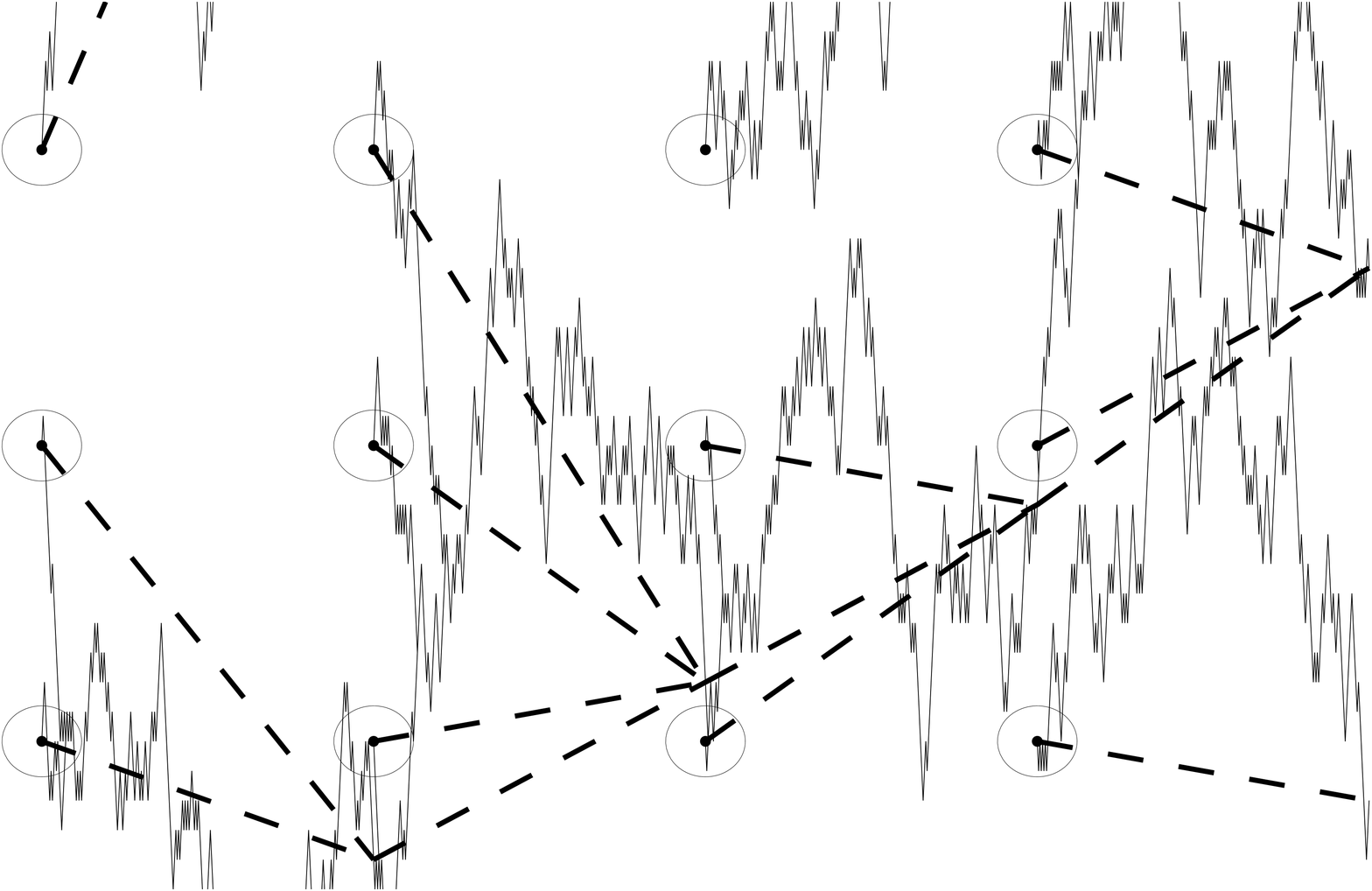}}
\put(8.7,0){\makebox(0,0){(b)}}
\end{picture}
\caption[]{\label{fig2}\small
A two-dimensional array of independent random signs can govern a
system of coalescing random walks (a) converging to a system of
coalescing Brownian motions (b).}
\end{center}
\end{figure}

Independent random signs ($ \pm 1 $, equiprobable) can govern various
interesting discrete models that have continuous scaling (that is,
mesh refinement) limits, see Figs~\ref{fig1}, \ref{fig2}.

Each real-valued function $ \phi $ of a finite collection of random
signs $ \tau_1, \dots, \tau_n $ can be written as a polynomial,
\[
\phi (\tau_1, \dots, \tau_n) = \hat \phi () + \sum_k \hat \phi (k)
\tau_k + \sum_{k<l} \hat \phi (k,l) \tau_k \tau_l + \dots + \hat \phi
(1,\dots,n) \tau_1 \dots \tau_n \> ;
\]
$ \hat \phi : 2^{\{1,\dots,n\}} \to \R $ is called the Fourier-Walsh
transform of $ \phi : \{-1,+1\}^n \to \R $. A continuous-time
counterpart, It\^o's decomposition into multiple stochastic integrals,
\[
\phi \( B(\cdot) \) = \hat \phi_0 + \int \hat \phi_1(t) \, dB(t) +
\iint\limits_{t_1<t_2} \hat \phi_2 (t_1,t_2) \, dB(t_1) dB(t_2) +
\dots \> ,
\]
works only for the classical case of Brownian motion (white noise). In
general, continuous time theory is more complicated; $ \hat \phi_n $
with finite $ n $ are not enough, and we need so-called spectral
decomposition of a noise, introduced in \cite{unitary}.

The main result of the present work (Theorem \ref{th4.2}) states,
roughly speaking, that the spectral decomposition of a noise is the
scaling limit of the Fourier-Walsh expansion for the corresponding
discrete model. An exact formulation is given after some definitions
that set a framework for the concept of scaling limit. Applications to
specific models (see \cite{coefs}, \cite{W}) will be published
separately.

%% file: first.tex
\section{Beyond topological semigroups}

Brownian motions in Lie groups are a classical topic. Brownian motions
in some infinite-dimensional topological groups arise naturally from
stochastic flows with smooth coefficients. However, we are mostly
interested in essentially non-smooth (coalescing, splitting, etc.)
stochastic flows. Corresponding maps (say, $ \R^1 \to \R^1 $) are far
from being invertible. Typically, such a map is a piecewise constant
function, defined (and continuous) everywhere except for some discrete
set. Accordingly, the composition $ g \circ f $ of such functions is
defined not for all pairs $ (f,g) $, but only for almost all pairs
(w.r.t.\ some relevant measures). The functions fail to form a
topological semigroup.

\begin{Def}\label{under}
An \emph{undergroup} $ G $ is a set equipped with a metric $ \rho $
and a map $ (f,g) \mapsto fg $ from a subset of $ G \times G $ to $ G
$, satisfying the following conditions.

(a) ``Unity'': there is $ e \in G $ such that for every $ f \in G $
both $ ef $ and $ fe $ are defined, and $ ef = f = fe $.

(b) ``Associativity'': let $ f,g,h \in G $ be such that $ fg $ and $
gh $ are defined; then $ (fg)h $ and $ f(gh) $ are defined, and $
(fg)h = f(gh) $.

(c) ``Continuity'': let $ f,g, f_1,g_1, f_2,g_2, \dots \in G $ be such
that $ fg $, $ f_1 g_1 $, $ f_2 g_2 $, \dots\ are defined; if $ \rho
(f_n,f) \to 0 $ and $ \rho (g_n,g) \to 0 $, then $ \rho ( f_n g_n, fg
) \to 0 $.

(d) Let $ f,g \in G $ be such that $ fg $ is defined; then $ \rho
(f,fg) \le \rho (e,g) $ and $ \rho (g,fg) \le \rho (e,f) $ (here $ e $
is the unity stipulated by (a); it is evidently unique).

(e) $ \rho (f,g) \le 1 $ for all $ f,g \in G $.
\end{Def}

\begin{example}\label{example1}
Let $ H_n $ be the set of all $ n $-point subsets of $ (0,1) $, and $
G = \cup_{n=0}^\infty H_n \times H_{n+1} $. An element $ f = (A,B) =
\( \{ a_1,\dots,a_n \}, \{ b_0,b_1,\dots,b_n \} \) $ of $ G $ may be
treated as a function $ f : (0,1) \setminus \{a_1,\dots,a_n\} \to
(0,1) $ defined by $ f(x) = b_k $ for all $ x \in (a_k,a_{k+1}) $,
where $ k = 0,1,\dots,n $, and $ a_0 = 0 $, $ a_{n+1} = 1 $; see
Fig.~\ref{fig3}a.

\begin{figure}[b]
\begin{center}
\setlength{\unitlength}{1cm}
\begin{picture}(13,2.6)
\put(-0.3,0.4){\includegraphics{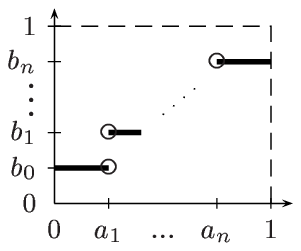}}
\put(1.5,0){\makebox(0,0){(a)}}
\put(3.7,0.2){\includegraphics{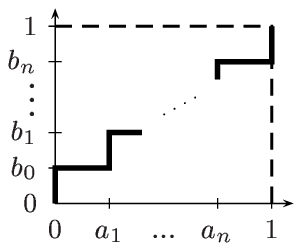}}
\put(5.5,0){\makebox(0,0){(b)}}
\put(9.1,0.25){\includegraphics{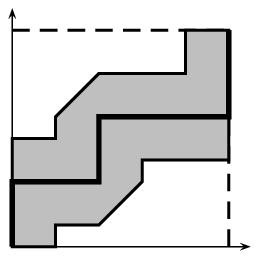}}
\put(10.5,0){\makebox(0,0){(c)}}
\end{picture}
\caption[]{\label{fig3}\small
An element of $ G $ as a function (a) or a (many-to-many)
binary relation (b); a \mbox{$ \rho $-neighborhood} (c).}
\end{center}
\end{figure}

Composition $ fg $ of $ f=(A,B) $ and $ g=(C,D) $ is defined if and
only if $ B \cap C = \emptyset $, and corresponds to the usual
composition of functions, $ x \mapsto g(f(x)) $. It is easy to see
that $ fg $ is of the form $ (A',D') $ for some $ A' \subset A $ and $
D' \subset D $.

Alternatively, an element $ f = (A,B) = \( \{ a_1,\dots,a_n \}, \{
b_0,b_1,\dots,b_n \} \) $ of $ G $ may be treated as a subset $ \{0\}
\times [0,b_0] \cup [0,a_1]\times \{b_0\} \cup \{a_1\}\times[b_0,b_1]
\cup \dots \cup [a_n,1]\times\{b_n\} \cup \{1\}\times[b_n,1] $ of $
[0,1] \times [0,1] $, see Fig.~\ref{fig3}b.  The subset belongs to the
metric space of all closed subsets of $ [0,1] \times [0,1] $ with
Hausdorff metric corresponding to $ l_1 $ metric $ \rho_1 \( (a,b),
(c,d) \) = |a-c| + |b-d| $ on the square. The Hausdorff metric induces
the metric $ \rho $ on $ G $, see Fig.~\ref{fig3}c.  The composition $
fg $ in $ G $ may be interpreted as composition of (many-to-many)
binary relations.

We add to $ G $ a unit $ e $, represented by $ f(x) = x $, or by the
diagonal of the square.

It can be shown that such $ G $ is an undergroup.
\end{example}

Unless stated otherwise, a measure on $ G $ is assumed to be a
probability measure concentrated on a countable union of $ \rho
$-compact subsets.

\begin{Def}\label{conv}
Let $ G $ be an undergroup, and $ \mu, \nu $ measures on $ G $. If the
composition $ fg $ is defined for $ \mu \otimes \nu $-almost all pairs
$ (f,g) $, then the \emph{convolution} $ \mu * \nu $ is defined as the
image of $ \mu \otimes \nu $ under $ (f,g) \mapsto fg $. Otherwise,
the convolution is undefined.
\end{Def}

The measure $ \mu * \nu $ (if defined) is concentrated on a countable
union of compact sets due to Lusin continuity of the map $ (G,\mu)
\times (G,\nu) \ni (f,g) \mapsto fg \in G $.

Condition \ref{under}(e), if violated, can be forced by replacing $
\rho (f,g) $ with $ \rho_1 (f,g) = \min ( 1, \rho(f,g) ) $. Due to
\ref{under}(e), weak convergence of measures on $ G $ may be metrized
by the transportation (Kantorovich-Rubinstein) metric
\[
\rhoKR (\mu,\nu) = \sup_\phi \bigg| \int \phi \, d\mu - \int \phi \,
d\nu \bigg| \, ,
\]
where the supremum is taken over all $ \phi : G \to \R $ satisfying
Lipshitz condition: $ | \phi(f) - \phi(g) | \le \rho(f,g) $ for all $
f,g \in G $.

\begin{lemma}\label{lemma1.4}
Let $ G $ be an undergroup. Then the set $ M(G) $ of all measures on $
G $, equipped with the transportation metric and the convolution, is
an undergroup.
\end{lemma}

\textsc{Proof.}
(a) The unit mass $ \de_e $ at $ e $ is the unity of $ M(G) $.

(b) Let $ \la * \mu $ and $ \mu * \nu $ be defined, then $ fg $ and $
gh $ are defined for $ \la\otimes\mu\otimes\nu $-almost all $ (f,g,h)
$. Therefore, $ fgh = (fg)h = f(gh) $ is defined; its distribution is
$ \la*\mu*\nu = (\la*\mu)*\nu = \la*(\mu*\nu) $.

(c) Let $ \rhoKR (\mu_n,\mu) \to 0 $, then we can construct $ G
$-valued random variables $ X, X_n $ on a probability space such that
$ X \sim \mu $, $ X_n \sim \mu_n $, and $ X_n \to X $ almost sure. The
same for $ \nu_n $, $ \nu $ and $ Y_n $, $ Y $ on another probability
space. Let $ \mu*\nu $ and $ \mu_n * \nu_n $ be defined. Then, on the
product of the two probability spaces, $ XY $ and all $ X_n Y_n $ are
defined almost sure, and $ X_n Y_n \to XY $ almost sure, since 1.1(c)
holds for $ G $. However, $ XY \sim \mu*\nu $ and $ X_n Y_n \sim \mu_n
* \nu_n $. Therefore $ \rhoKR ( \mu_n * \nu_n, \mu * \nu ) \to 0 $.

(d) Let $ \mu*\nu $ be defined, then $ fg $ is defined for $
\mu\otimes\nu $-almost all $ (f,g) $, and $ \rho (f,fg) \le \rho (e,g)
$. Let $ \phi : G \to \R $ satisfy $ | \phi(f) - \phi(g) | \le \rho
(f,g) $ for all $ f,g $. We have $ | \phi(f) - \phi(fg) | \le \rho
(f,fg) \le \rho(e,g) $, therefore
\[\begin{split}
& \rhoKR ( \mu, \mu*\nu ) = \sup_\phi \bigg| \iint ( \phi(f) -
  \phi(fg) ) \, d\mu(f) \, d\nu(g) \bigg| \le \\
& \le \iint \rho (e,g) \, d\mu(f) \, d\nu(g) = \int \rho
  (e,g) \, d\nu(g) \le \rhoKR ( \de_e, \nu ) \, ,
\end{split}\]
since the function $ \tilde\phi (g) = \rho (e,g) $ satisfies $ |
\tilde\phi (g_1) - \tilde\phi (g_2) | \le \rho (g_1,g_2) $, and so, $
\int \rho (e,g) \, d\nu(g) = | \int \tilde\phi (g) \, d\nu - \int
\tilde\phi (g) \, d\de_e | \le \rhoKR ( \nu, \de_e ) $. Similarly, $
\rhoKR ( \nu, \mu*\nu ) \le \rhoKR ( \de_e, \mu ) $.

(e) Let $ \phi $ satisfy $ | \phi (f) - \phi (g) | \le \rho (f,g) $,
then $ | \phi(f) - \phi(g) | \le 1 $, thus all values of $ \phi $ lie
inside some $ [A,A+1] $. Therefore, $ \int \phi \, d\mu $ and $ \int
\phi \, d\nu $ also belong to $ [A,A+1] $, and so, $ \rhoKR (\mu,\nu)
= \sup_\phi | \int \phi \, d\mu - \int \phi \, d\nu | \le 1 $.\qed

%% file: second.tex
\section{Convolution semigroups and independence (continuous time)}

A convolution semigroup, defined below, is just a one-parameter
semigroup in $ M(G) $. A one-parameter semigroup in $ G $ could be
defined, but is of little interest, since the undergroup of Example
\ref{example1} (unlike a Lie group) contains no nontrivial
one-parameter semigroup. Nevertheless, it contains nontrivial (and
interesting) convolution semigroups.

\begin{Def}\label{def1.5}
Let $ G $ be an undergroup. A \emph{convolution semigroup} in $ G $ is
a family $ \(\mu_t\)_{t\in[0,\infty)} $ of measures $ \mu_t \in M(G) $
such that

(a) $ \mu_0 $ is the unit of $ M(G) $;

(b) $ \rhoKR ( \mu_t, \mu_0 ) \to 0 $ for $ t \to 0 $;

(c) for every $ s,t \in [0,\infty) $, the convolution $ \mu_s * \mu_t
$ is defined, and $ \mu_s * \mu_t = \mu_{s+t} $.
\end{Def}

Note that \ref{def1.5}(b) combined with \ref{lemma1.4}(d) ensures
continuity of $ \mu_t $ in $ t $.

\begin{lemma}\label{lemma2.1}
Let $ \(\mu_t\) $ be a convolution semigroup in an undergroup $ G
$. Then there exists a two-parameter family $ \( X_{s,t} \)_{0\le s
\le t < \infty} $ of $ G $-valued random variables on some probability
space ($ (0,1) $ with Lebesgue measure can be used) such that for
every $ s,t $ satisfying $ 0 \le s \le t < \infty $,

(a) $ \mu_{t-s} $ is the distribution of $ X_{s,t} $,

(b) $ X_{r,s} X_{s,t} = X_{r,t} $ almost sure,\footnote{%
The exceptional set may depend on $ r,s,t $.}

\noindent and for every $ n $ and $ t_1,\dots,t_n $ satisfying $ 0 \le
t_1 \le \ldots \le t_n < \infty $,

(c) $ X_{0,t_1} $, $ X_{t_1,t_2} $, \dots, $ X_{t_{n-1},t_n} $ are
independent.
\end{lemma}

The proof is left to the reader.

The family $ (X_{s,t}) $ is a $ G $-valued counterpart of the
classical ``process with independent increments''. Unlike the
classical case, $ X_{0,s} X_{s,t} = X_{0,t} $ cannot be written as $
X_{s,t} = (X_{0,s})^{-1} X_{0,t} $. These $ X_{s,t} $ are independent,
but they are not increments. Such a family $ (X_{s,t}) $ may be called
a $ G $-valued \emph{independent process.} It determines \sif
s\footnote{%
Every \sif\ is assumed to contain all negligible sets.}
$ \F_{s,t} $ for $ 0 \le s \le t < \infty $; namely, $ \F_{s,t} $ is
generated by $ G $-valued random variables $ X_{u,v} $ for all $ u,v $
satisfying $ s \le u \le v \le t $. Clearly,
\[
\F_{r,s} \otimes \F_{s,t} = \F_{r,t}
\]
whenever $ 0 \le r \le s \le t < \infty $; that is, $ \F_{r,s} $ and $
\F_{s,t} $ are independent and, taking together, they generate $
\F_{r,t} $.

The construction may be extended from $ [0,\infty) $ to $
(-\infty,+\infty) $. To this end, introduce $ X_{s,t} $ for $ -\infty
< s \le t \le 0 $ satisfying (negative-time counterparts of)
\ref{lemma2.1}(a,b,c) and independent of those introduced before (for
positive time). Let $ X_{s,t} = X_{s,0} X_{0,t} $ whenever $ -\infty <
s < 0 < t < +\infty $. The extended construction is invariant under
time shifts:
\[
X_{s+u,t+u} = X_{s,t} \circ T_u
\]
where $ \( T_u \)_{u\in\R} $ is a one-parameter group of measure
preserving transformations of our probability space.\footnote{%
Choose a Lebesgue space (say, $ (0,1) $) for the probability space,
and assume that $ \F_{-\infty,+\infty} $ contains all measurable sets
(otherwise a quotient space should be used), then existence of $ (T_u)
$ follows from the fact that the joint distribution of $
X_{t_0+u,t_1+u} $, \dots, $ X_{t_{n-1}+u,t_n+u} $ does not depend on $
u $.}
We get a noise, as defined below, following \cite[Def.~1.1]{unitary}.

\begin{Def}\label{noise}
A \emph{noise} consists of a probability space $ (\Om,\F,P) $, a
one-parameter group $ \(T_t\)_{t\in\R} $ of measure preserving
transformations $ T_t : \Om \to \Om $, and a two-parameter family $
(\F_{s,t}) $ of sub-\sif s $ \F_{s,t} \subset \F $ for $ -\infty < s
\le t < \infty $, such that for all $ r,s,t\in \R $

(a) $ T_t $ sends $ \F_{r,s} $ onto $ \F_{r+t,s+t} \, $ ($ r \le
s $);

(b) $ \F_{r,s} $ and $ \F_{s,t} $ are independent ($ r \le s \le
t $);

(c) $ \F_{r,s} $ and $ \F_{s,t} $, taken together, generate $ \F_{r,t}
\, $ ($ r \le s \le t $).
\end{Def}

%% file: third.tex
\section{Discrete time: convolution, independence, spectrum}

Notions discussed in Section 2 have evident discrete-time
counterparts. We start from an undergroup $ G $ and a measure $ \mu_1
\in M(G) $ such that $ \mu_1 * \mu_1 $ is defined; we introduce
measures $ \mu_t = \underbrace{ \mu_1 * \ldots * \mu_1 }_{t} $ for $ t
= 1,2,\ldots $ ($ \mu_0 $ being the unit of $ M(G) $) and independent
$ G $-valued random variables $ X_{t,t+1} $ for $ t \in \Z $,
distributed $ \mu_1 $. (Thus, $ (G,\mu_1)^\Z $ may
be chosen for the probability space.) We define $ X_{s,s} = e $ (the
unit of $ G $), and $ X_{s,t} = X_{s,s+1} \ldots X_{t-1,t} $ for $ s,t
\in \Z $, $ s < t $. The \sif\ $ \F_{s,t} $ is generated by $
X_{s,s+1}, \ldots, X_{t-1,t} $. Time shifts $ T_u = (T_1)^u $ satisfy
$ X_{s,t} \circ T_u = X_{s+u,t+u} $.

The simplest nontrivial case appears when $ \mu_1 $ is concentrated at
two equiprobable points (atoms) $ f_1, f_2 $. In that case, introduce
independent random signs $ \tau_t $, $ t \in \Z $, as follows:
\[
\tau_t = \begin{cases}
  -1, & \text{if $ X_{t-1,t} = f_1 $},\\
  +1, & \text{if $ X_{t-1,t} = f_2 $}.
\end{cases}
\]
Then $ X_{0,t} = X_{0,1} \dots X_{t-1,t} $ may be treated as a $ G
$-valued function of $ \tau_1, \ldots, \tau_t $. Given a function $
\psi : G \to \R $ and a natural $ n $, we get a real-valued function $
\phi $ of $ n $ random signs,
\[
\psi (X_{0,n}) = \phi ( \tau_1, \dots, \tau_n ) \, ,
\]
and its Fourier-Walsh transform $ \hat\phi : 2^{\{1,\dots,n\}} \to \R $,
\begin{equation}\label{3.1}
\phi (\tau_1,\ldots,\tau_n) = \hat\phi () + \sum_k \hat\phi (k) \tau_k
+ \sum_{k<l} \hat\phi (k,l) \tau_k \tau_l + \ldots + \hat\phi
(1,\ldots,n) \tau_1 \ldots \tau_n \, .
\end{equation}

Assume for convenience that $ \E | \psi (X_{0,n}) |^2 = 1 $, then
\[
|\hat\phi()|^2 + \sum_k |\hat\phi(k)|^2 + \sum_{k<l} |\hat\phi(k,l)|^2
 + \ldots + |\hat\phi(1,\ldots,n)|^2 = 1 \, ;
\]
the summands may be treated as probabilities that describe a random
subset $ \Sra $ of $ [0,n] $.\footnote{%
 Usually, all random variables are functions on a single probability
 space, they have a joint distribution. However, $ \Sra $ has no joint
 distribution with $ \tau_t $ (or $ X_{s,t} $). It may be thought of as
 a manifestation of quantum complementarity. Functions $ \phi $ and $
 \hat\phi $ may be treated as coefficients of a quantum state vector in
 two different orthonormal bases. Moreover, if $ \phi $ describes a
 spin state of an $ n $-electron system, then $ \hat\phi $ describes
 the same state in a coordinate system rotated by $ 90^\circ $ in the
 three-dimensional space, and $ \Sra $ is the outcome of a quantum
 measurement, governed by $ |\hat\phi(\dots)|^2 $, incompatible with
 another measurement, governed by $ |\phi(\dots)|^2 $.}
Namely,
\[\begin{split}
& |\hat\phi()|^2 = \Pr{ \Sra = \emptyset } \, ; \\
& |\hat\phi(k)|^2 = \Pr{ \Sra = [k-1,k] } \, ; \\
& |\hat\phi(k,l)|^2 = \Pr{ \Sra = [k-1,k] \cup [l-1,l] } \, ; \\
& \dots \\
& |\hat\phi(1,\dots,n)|^2 = \Pr { \Sra = [0,n] } \, .
\end{split}\]
For now, these probabilities are defined only for the special case of
a two-atomic $ \mu_1 $. A generalization is suggested by the formula
\begin{equation}\label{3.2}
\Pr { \Sra \subset [u,v] } = \E \,\big|\, \cE{ \psi(X_{0,n}) }{
  \F_{u,v} } \,\big|^2 \, ;
\end{equation}
here the function $ \psi (X_{0,n}) = \phi (\tau_1,\dots,\tau_n) $ is
averaged over $ \tau_1, \dots, \tau_u $ and $ \tau_{v+1}, \dots,
\tau_n $, giving a function of $ \tau_{u+1}, \dots, \tau_v $; the
latter function is squared and then averaged over $ \tau_{u+1}, \dots,
\tau_v $. The \sif\ $ \F_{u,v} $ is generated by $ X_{t-1,t} $ for all
$ t $ such that $ [t-1,t] \subset [u,v] $.

The proof of \eqref{3.2} is easy. Write \eqref{3.1} in the form
\[
\phi (\tau_1,\dots,\tau_n) = \sum_{S\subset\{1,\dots,n\}} \hat\phi(S)
\tau_S \, , \quad \text{where } \tau_S = \prod_{k\in S} \tau_k \, ,
\]
and note that $ \cE{ \tau_S }{ \F_{u,v} } = \tau_S $ if $ S \subset \{
u+1, \dots, v \} $, othervise $ \cE{ \tau_S }{ \F_{u,v} } $ $ = 0 $;
we have
\[\begin{split}
& \cE{ \phi (\tau_1,\dots,\tau_n) }{ \F_{u,v} } =
  \sum_{S\subset\{u+1,\dots,v\}} \hat\phi (S) \tau_S \, , \\
& \E \,|\dots|^2 = \sum_{S\subset\{u+1,\dots,v\}} | \hat\phi(S) |^2 = \Pr
  { \Sra \subset [u,v] } \, ,
\end{split}\]
which proves \eqref{3.2}.

Equality \eqref{3.2} describes the joint distribution of the minimal
and the maximal elements of $ \Sra $, not the whole $ \Sra $. However,
we may introduce the \sif\ $ \F_{u_1,v_1,u_2,v_2} $ (for $ 0 \le u_1 <
v_1 < u_2 < v_2 \le n $) generated by $ X_{t-1,t} $ for all $ t $ such
that $ [t-1,t] \subset [u_1,v_1] \cup [u_2,v_2] $. Then
\[
\Pr { \Sra \subset [u_1,v_1] \cup [u_2,v_2] } = \E \,\big|\, \cE{
\psi(X_{0,n}) }{ \F_{u_1,v_1,u_2,v_2} } \big|^2 \, .
\]
Generalization for $ E = [u_1,v_1] \cup \dots \cup [u_k, v_k] $ is
now straightforward:
\begin{equation}\label{3.3}
\Pr { \Sra \subset E } = \E \,\big|\, \cE{ \psi(X_{0,n}) }{ \F_E }
\big|^2 \, .
\end{equation}
That is enough for describing the distribution of $ \Sra $ via an
inclusion-exclusion formula
\[
\Pr { \Sra = S } = \sum_{E\subset S} (-1)^{|S|-|E|} \, \Pr { \Sra
\subset E } \, .
\]
Here $ S $ and $ E $ run over sets of the form $ [u_1,v_1] \cup \dots
\cup [u_k,v_k] $, $ k \in \{ 0, 1, \dots, n \} $, $ u_i,v_i \in \{
0,1,\dots,n \} $, $ 0 \le u_1 < v_1
< u_2 < \dots <v_k \le n $, and $ |S| $ is Lebesgue measure of $ S $.

The distribution of $ \Sra $, defined by \eqref{3.3} for an arbitrary
$ \mu_1 $, will be called the \emph{spectral measure} of the random
variable $ \psi (X_{0,n}) $.

%% file: fourth.tex
\section{Scaling limit}

Let $ \( \mu_t \)_{t\in[0,\infty)} $ be a (continuous-time)
convolution semigroup (as defined by 2.1) in an undergroup $ G $.
On the other hand, let $ \mu^{(n)} \in M(G) $ be given for $ n =
1,2,\dots $ such that $ \mu^{(n)} * \mu^{(n)} $ is defined. We
introduce measures $ \mu^{(n)}_{k/2^n} = \underbrace{ \mu^{(n)} *
\dots * \mu^{(n)} }_k $ that form discrete-time convolution semigroups
$ \( \mu^{(n)}_t \)_{t\in2^{-n}\Z_+} $ similarly to Sect.~3, but the $
n $-th semigroup has its time pitch $ 2^{-n} $. Assume that
\begin{equation}\label{4.1}
\rhoKR \( \mu^{(n)}_t, \mu_t \) \to 0 \quad \mbox{for $ n\to\infty $}
\end{equation}
for every binary-rational (that is, of the form $ k \cdot 2^{-l} $)
number $ t \ge 0 $.

Let a function $ \psi : G \to \R $ be continuous $ \mu_1 $-almost
everywhere and bounded, and $ \int_G | \psi(f) |^2 \, \mu_1(df) = 1
$. Introduce $ c_n = \int_G | \psi(f) |^2 \, \mu^{(n)}_1(df) $, then $
c_n \to 1 $. For each $ n $ we repeat the construction of Sect.~3, but
with time pitch (mesh) $ 2^{-n} $ rather than $ 1 $. We get $ G $-valued
random variables\footnote{%
Underlying probability spaces may depend on $ n $.}
$ X^{(n)}_{s,t} $ for $ s,t \in 2^{-n}\Z $, $ s < t $, and $ \E \,|
\psi ( X^{(n)}_{0,1} ) |^2 = c_n $. The spectral measure $ \nu_n $ of
the random variable $ c_n^{-1/2} \psi ( X^{(n)}_{0,1} ) $ is the
distribution of a random set $ \Sra_n \subset [0,1] $ that respects
the time pitch $ 2^{-n} $.

If there is a limit of $ \nu_n $ for $ n \to \infty $, it should be
the distribution of a random subset $ \Sra $ of $ [0,1] $, so that $
\Sra_n \to \Sra $ in distribution. To this end, however, we need an
appropriate metric for subsets. Hausdorff metric
\[
\rho_\mathrm{H} ( S_1, S_2 ) = \inf \{ \eps>0 : S_1 \subset
\(S_2\)_{+\eps}, S_2 \subset \(S_1\)_{+\eps} \}
\]
will be used; here $ S_1, S_2 \subset
[0,1] $ are closed sets, and $ \( S \)_{+\eps} $ is the $ \eps
$-neighbor\-hood of $ S $.

\begin{theorem}\label{th4.2}
(a) The sequence $ \( \Sra_n \)_{n=1,2,\dots} $ converges in
distribution in the
Hausdorff space (of all closed subsets of $ [0,1] $, with Hausdorff
metric) to some $ \Sra $.

(b) Almost surely, $ \Sra $ is a closed subset of $ [0,1] $ having
zero Lebesgue measure.

(c) $ \Pr { t \in \Sra } = 0 $ for every $ t \in [0,1] $.

(d) For every finite union $ E \subset [0,1] $ of intervals,
\[
\Pr { \Sra \subset E } = \E \,|\, \cE{ \psi(X_{0,1}) }{ \F_E } |^2 \,
.
\]
\end{theorem}

\textsc{Proof.}
(a,d): The Hausdorff space is compact, therefore the space of all
probability distributions on that space is also compact. It suffices
to prove that the sequence $ (\mu_n) $ has only one limit
point. However, a distribution on the Hausdorff space is uniquely
determined by probabilities of the form $ \Pr { \Sra \subset E } $
where $ E $ is a finite union of intervals with binary-rational
endpoints. Therefore it suffices to prove convergence of $ c_n \Pr {
\Sra_n \subset E } = \E \,|\, \cE{ \psi (X_{0,1}^{(n)}) }{ \F_E } |^2
$ for $ n \to \infty $.

Start with a special case: $ E = [0,t] $ for some binary-rational $ t
\in (0,1) $. Due to \eqref{4.1}, we can construct $ G $-valued random
variables $ X_{0,t} $ and $ X_{0,t}^{(n)} $ on some probability space
such that $ X_{0,t} \sim \mu_t $, $ X_{0,t}^{(n)} \sim \mu_t^{(n)} $,
and $ X_{0,t}^{(n)} \to X_{0,t} $ almost sure. Similarly, we construct
$ X_{t,1} $ and $ X_{t,1}^{(n)} $ on another probability space such
that $ X_{t,1} \sim \mu_{1-t} $, $ X_{t,1}^{(n)} \sim \mu_{1-t}^{(n)}
$, and $ X_{t,1}^{(n)} \to X_{t,1} $ almost sure.  Denote by $
\F_{0,t} $ the \sif\ of the first probability space, by $ \F_{t,1} $
the \sif\ of the second, and by $ \F_{0,1} = \F_{0,t} \otimes \F_{t,1}
$ the \sif\ of their product. On the product space, introduce $ G
$-valued random variables $ X_{0,1}^{(n)} = X_{0,t}^{(n)}
X_{t,1}^{(n)} \sim \mu_1^{(n)} $ and $ X_{0,1} = X_{0,t} X_{t,1} \sim
\mu_1 $. We have $ X_{0,1}^{(n)} \to X_{0,1} $ almost sure by
\ref{under}(c). However, $ \psi $ is continuous $ \mu_1 $-almost
everywhere and bounded, thus, $ \psi (X_{0,1}^{(n)}) \to \psi
(X_{0,1}) $ almost sure, and $ \cE{ \psi (X_{0,1}^{(n)}) }{ \F_{0,t} }
\to \cE{ \psi (X_{0,1}) }{ \F_{0,t} } $ almost sure. Therefore $ \E
\,|\, \cE{ \psi (X_{0,1}^{(n)}) }{ \F_{0,t} } |^2 \to \E \,|\, \cE{
\psi (X_{0,1}) }{ \F_{0,t} } |^2 $. It means convergence of $ \Pr {
\Sra_n \subset [0,t] } $.

Instead of the general case, consider another special case $ E =
[u_1,v_1] \cup [u_2,1] $, since further generalization is
straightforward. We get convergence almost sure, $ X^{(n)}_{0,u_1}
\to X_{0,u_1} $ on some probability space, $ X^{(n)}_{u_1,v_1}
\to X_{u_1,v_1} $ on another, $ X^{(n)}_{v_1,u_2} \to X_{v_1,u_2} $ on
a third, and $ X^{(n)}_{u_2,1} \to X_{u_2,1} $ on a fourth probability
space. We multiply the four spaces and consider the \sif\ $ \F_E =
\F_{u_1,v_1} \otimes \F_{u_2,1} $ generated by second and fourth
spaces. Convergence is established as before, which completes the
proof of (a) and (d).

(c): Applying (d) to $ E_\eps = [0,1] \setminus (t-\eps, t+\eps) $ we
reduce (c) to the statement that the \sif\ $ \F_{0,1} $ is generated
by the union of \sif s $ \F_{E_\eps} $ for all $ \eps > 0 $. Note that
$ \rho ( X_{0,t-\eps}, X_{0,t} ) \to 0 $ in probability by
\ref{def1.5}(b). Choose $ \eps_n \to 0 $ such that $ \rho (
X_{0,t-\eps_n}, X_{0,t} ) \to 0 $ almost sure and also $ \rho (
X_{t+\eps_n,1}, X_{t,1} ) \to 0 $ almost sure. It follows by
\ref{under}(c) that $ X_{0,t-\eps_n} X_{t+\eps_n,1} \to X_{0,t}
X_{t,1} = X_{0,1} $, therefore $ X_{0,1} $ is measurable w.r.t.\ the
\sif\ generated by all $ \F_{E_\eps} $. The same holds for $ X_{r,s}
$, $ 0 \le r \le s \le 1 $, which gives the whole $ \F_{0,1} $.

(b): Follows immediately from (c).\qed

\begin{note}
The distribution of $ \Sra $ is a special case of
\cite[(2.8)]{unitary}. It is uniquely determined (via \ref{th4.2}(d))
by the convolution semigroup $ (\mu_t) $ and the function $ \psi $,
irrespective of any discrete approximations. For the classical case of
a Brownian motion in $ \R^n $ or a Lie group, the set $ \Sra $ is
finite (almost sure). In general, finiteness of $ \Sra $ is necessary
and sufficient for the noise to be classical (white) up to
isomorphism, see \cite[2.14]{unitary}.
\end{note}

\begin{note}
Properties \ref{th4.2}(b,c) are known for all noises, see
\cite[2.3(a)]{unitary}.
\end{note}

\begin{note}
For the first nonclassical noise \cite[Sect.~5]{TV} it is still
unknown, whether it corresponds to some convolution semigroup in some
undergroup, or not. When writing \cite{TV}, Anatoly Vershik asked me
repeatedly about something like a semigroup structure.
\end{note}

%% file: main.bbl
\begin{thebibliography}{99}

\bibitem{coefs} B.~Tsirelson, ``Fourier-Walsh coefficients for a
coalescing flow (discrete time)'', math.PR/9903068.

\bibitem{unitary} B.~Tsirelson, ``Unitary Brownian motions are
linearizable'', math.PR/9806112.

\bibitem{TV} B.S.~Tsirelson, A.M.~Vershik,
``Examples of nonlinear continuous tensor products of measure spaces
and non-Fock factorizations'',
Reviews in Mathematical Physics \textbf{10}:1 (1998), 81--145.

\bibitem{W} J.~Warren, ``The noise made by a Poisson snake'',
Manuscript, Univ.\ de Pierre et Marie Curie, Paris, Nov.\ 1998.

\end{thebibliography}
